\documentstyle[12pt]{article}
\begin{document}
\setlength{\textwidth}{6in}
\setlength{\textheight}{8in}
\setlength{\topmargin}{0.3in}
\setlength{\headheight}{0in}
\setlength{\headsep}{0in}
\setlength{\parskip}{0pt}
\input{mssymb}

\newcounter{probnumt}
\partopsep 0pt
\settowidth{\leftmargin}{1.}\addtolength{\leftmargin}{\labelsep}
\newenvironment{probt}{\begin{quote}\begin{enumerate} 
\setcounter{enumi}
{\value{probnumt}}}%
{  \setcounter{probnumt}
{\value{enumi}}\end{enumerate}\end{quote}}

\newcommand{\RR}{{\bf R}}
\newcommand{\CC}{{\bf C}}
\newcommand{\NN}{{\bf N}}
\renewcommand{\limsup}{\overline {\lim}\,}
\renewcommand{\liminf}{\underline {\lim}\,}
\newtheorem{thm}{Theorem}
\newtheorem{lem}[thm]{Lemma}
\newtheorem{cor}[thm]{Corollary}
\newtheorem {ex}[thm]{Example}
\newtheorem{pro}{Question}
\newenvironment{rem}{\medskip\par\noindent {\bf Remark 1.}}

\newenvironment{proof}{\medskip\par\noindent{\bf Proof.}}{\hfill
$\Box$ \medskip \par}

\setcounter{page}{1}

\title{Extremal properties of contraction semigroups on $c_o$}

\author { Pei-Kee Lin\\
Department of Mathematics \\
University of Memphis \\
Memphis TN, 38152\\
USA}

\maketitle
\footnotetext[3]{partial results of this paper obtained when the author
attained the International Conference on Convexity at the
University of Marne-La-Vall\'ee.  He would like to thank the kind
hosptitality offered to him.
He also
likes to think Professor Goldstein and Professor 
Jamison for their valuable
suggestions. }

\begin{abstract}
{For any  complex Banach space $X$, let $J$ denote the duality
mapping of $X$.
For any unit vector  $x$  in $X$ and 
 any
($C_0$) contraction semigroup  $(T_t)_{t>0}$ on $X$,
Baillon and Guerre-Delabriere proved 
that if $X$ is a smooth reflexive Banach space and if there is
$x^* \in J(x)$ such that $|\langle T(t)
\, x,J(x)\rangle| \to 1 $ as $t \to \infty$, then there is 
a unit vector $y\in X$ which is an
eigenvector of the generator $A$ of $(T_t)_{t>0}$
associated with a purely imaginary eigenvalue. 
They asked whether this result is still true if $X$ is replaced
by $c_o$. In this article, we show the answer is negative.}
\end{abstract} \newpage

Let $X$ be a complex Banach space.
The {\em duality mapping} $J$ of $X$ is a (multivalued) function
from $X\setminus \{0\}$ into $X^*$ which is defined by
$$
J(x)=\{x^* \in X^*:\langle x,x^*\rangle=1=\|x^*\|\}.$$
Goldstein  \cite {Gb} proved the
following theorem.

\begin{thm} \label{A}  Let $A$ be a generator of a $C_0$ contraction
semigroup $T=\{e^{tA}:t \geq 0\}$ on a complex Banach  space $X$.
If $X$ is a Hilbert space, then for any unit vector $x \in X$,
\[\lim_{t \to \infty}|\langle T_t\, x, J(x) \rangle| = 1\]
implies that $A \, x=i \lambda x$ for some real number
$\lambda$. {\rm (}Note: since any  Hilbert space is a smooth
Banach space,  $J(x)$ is
singleton for $x \ne 0$.{\rm )} \end{thm}

He also showed that Theorem~\ref {A} is
not true if $X$ is replaced by an $L_\infty$-space.  In \cite{L}, 
the author proved that Theorem ~\ref {A} holds if and only if $X$
is strictly convex (also see \cite {BG}).
We shall note:  the counter-example in \cite {Gb}, the generator
of the $C_0$-semigroup has no eigenvector. On the other hand, the
counter-example in \cite {L}, the generator has an eigenvector
vector associated with a pure imaginary eigenvalue.  In \cite {BG},
Baillon and Guerre-Delabriere considered the following question.

\begin{pro} \label {Aa} Let $(T_t)_{t>0}$ be a $C_0$ contraction
semigroup on a complex Banach space $X$.
Suppose that $x$ is a unit vector in $X$ such that 
$$\lim_{t \to \infty}
|\langle T_t\,x, x^*\rangle|=1=\|x^*\|.$$
Find a necessary and sufficient condition so that there is a unit
vector $y \in X$ which is an eigenvector of $A$ associated with
a purely imaginary eigenvalue.
\end{pro}

They proved the following
theorems.

\begin{thm}\label {Ab} Let $(T_t)_{t >0}$ be a $C_0$
contraction semigroup  on
a complex Banach space $X$.  Let $x$ be a unit vector in $X$.
If there is $x^* \in J(x)$ such that
$$\lim_{t \to\infty} |\langle T_t \, x,x^*\rangle|=1,$$
then there is $y^* \in J(x)$ such that  $|\langle
T_t\, x, y^*\rangle|=1$ for all $t>0$.
\end{thm}

\begin{thm}\label {B}
Let $X$ be  a  reflexive smooth complex Banach space and let
$(T_t)_{t>0}$ be a semigroup of contractions on $X$.  If there
are a unit vector $x \in X$ and $x^* \in J(x)$ such that
$$
\lim_{t \to \infty}|\langle T_t \, x, x^*\rangle|=1,$$
then there are a unit vector $y \in X$ and $\omega\in \RR$ such
that for all $t>0$, $T_t \, y=e^{i\omega t}\, y$. Hence, $y$ is an
eigenvector of the generator of $(T_t)_{t>0}$ associated with an
eigenvalue $i\omega$. \end{thm}

They also showed Theorem ~\ref {B} is not true
if one replaces $X$  by  $\ell_1$ or $L_1$.  They asked
whether Theorem ~\ref {B} is still true if  $X$ is replaced by 
$c_o$ or any (non-smooth) reflexive Banach space.  
Recently, Ruess \cite {R} showed the answer is affirmative if $X$
is reflexive. Indeed, he proved the following strong theorem.

\begin{thm}  Assume that $(T_t)_{t>0}$ of a 
uniformly bounded  $C_0$-semigroup
on $X$.  If there exist $x \in X$ and $x^* \in X^*$ 
such that
$\{T_t \, x: t \in \Bbb R\}$ is weakly compact and
$\liminf_{t\to \infty}|\langle T_t\, x, x^*\rangle|>0$, then
there is a $y \in X\setminus \{0\}$ which is an eigenvector of
the generator of $(T_t)_{t>0}$ corresponding to a pure imagary
number.
\end{thm}

In this article, we consider the contraction $C_0$-semigroups on
$c_o$ and show the answer is negative.

\begin{ex}\label {C} {\rm Let $\{e_1,e_2,\cdots\}$ be the natural basis of
$c_o$.  $A:c_o \to c_o$ is defined by
$$
A(e_i)=\left \{ \begin{array} {ll}
  \sum_{k=2}^\infty \frac 1 k e_k \hspace{.3in} & \hbox {if $i=1$,}\\
   -\frac 1 i e_i & \hbox
  {otherwise.}\end{array}\right .$$
Then
$$
T_t(e_i)=e^{tA}(e_i)=\left \{ \begin{array} {ll}
  e_1 +\sum_{k=2}^\infty (1-e^{-\frac t k}) e_k\hspace{.3in} 
  & \hbox {if $i=1$,}\\
   e^{-\frac t i} e_i & \hbox
  {otherwise.}\end{array}\right .$$
So $(T_t)_{t>0}$ is a $C_0$ semigroup of contractions on $c_o$.
It is easy to see that 
$$\langle e^{tA} (e_1),e^*_1\rangle=1,$$
and $A$
does not have any eigenvector associated with any purely imaginary
eigenvalue (or 0). 
(Here $e^*_1 \in J(e_1)$ and $e^*_1$ is $e_1$, veiwed as a member
of $\ell_1={c_o}^*$.)
On the other hand, 
for any $k\geq 2$, $e_k$ is an eigenvector
of $A$ associated with eigenvalue $-\frac 1 k$.} \end{ex}

It is known that if $T$ is an isometry on $\ell_p$,
$1\leq p<\infty$, $p\ne 2$,
then the images $T\, x$ and $T\,y$
of any two disjoint elements $x,y$ (this is $|x|\wedge |y|=0$) 
are disjoint (\cite {Ro}
p. 416 Lemma 23). 
Let
$(T_t)_{t>0}$ be an isometric semigroup on $\ell_1$.
Since $\lim_{t\downarrow 0} T\, e_i =e_i$,
Fleming, Goldstein, and Jamison  \cite {FGJ,Gc} proved
there exists a sequence $\{\omega_n\}$ of $\Bbb R$ such that
$$
T_t(e_n)=e^{i\, \omega_n t} e_n.$$
Hence, there is no counter-examples of isometric semigroups on
$\ell_1$.
On the other hand, there is an isometry
$$T(\sum_{i=1}^\infty a_i \, e_i)= \frac {a_1+a_2} 2 \, e_1
+\sum _{i=1}^\infty a_i\, e_{i+1}$$
on $c_o$ 
which does not  preserve disjoint support.
It is natural to ask whether we can improve the Example
~\ref {C} to be $C_0$-isometric semigroup.  The
following theorem shows the answer is negative.

\begin{thm}  Let $(T_t)_{t>0}$ be a $C_0$ isometric semigroup on
$c_o$.  Let $\{e_k\}$ be the natural basis of $c_o$.  Then for
any $k \in \Bbb N$, there is $\omega_k \in \Bbb R$ such that
$$T_t\, e_k=e^{i\omega_kt} \, e_k.$$                \end{thm}

\begin{proof}  Let $\{e^*_k\}$ be the natural basis of
$\ell_1=(c_o)^*$ and let $k$ be any fixed natural number.  Since
$\lim_{t \downarrow 0} T_t\, e_k=e_k,$
there is $\delta_k>0$ such that if $0<t\leq \delta_k$, then
$$|\langle T_t\, e_k,e^*_j\rangle|<\frac 1 2 \hspace {.3in} \hbox
{for all $j \ne k$.}$$
But $T_t$ is an isometry.  This implies
$$|\langle T_t\, e_k,e_k^*\rangle|=1.$$
We claim that if $0<t\leq \delta_k$, then $\langle T_t\,
e_j,e^*_k\rangle=0$ for all $j \ne k$.  Suppose it is not true. 
There are $j \ne k$, and $0<t\leq \delta_k$, such that
$$\langle T_t\, e_j,e_k^*\rangle\ne 0.$$
So
$$\left\langle T_t\bigl(\overline{\langle T_t\, e_k,e_k^*\rangle}\,
e_k+\frac {
\overline{\langle T_t \, e_j,e^*_k\rangle}}{
|
\langle T_t \, e_j,e^*_k\rangle|}\, e_j\bigr),e_k^*\right\rangle =1+|
\langle T_t \, e_j,e^*_k\rangle|>1.$$
But
$$\left\|\overline{\langle T_t \, e_j,e^*_k\rangle} \, e_k+\frac{
\overline{\langle T_t \, e_j,e^*_k\rangle}}{|
\langle T_t \, e_j,e^*_k\rangle|}\, e_j\right\|=1.$$
This contradicts that $T_t$ is an isometry.  We proved our
claim.

We note: the claim proved that if $0<t<\delta_k$ and $x \in c_o$
such that $\langle x,e_k^*\rangle=0$, then $\langle T_t\,
x,e_k^*\rangle=0$.  Now, let $t$ be any positive real number. 
There exist integer $n \geq 0$ and $0 \leq s<\delta_k$ such that
$$t=n\, \delta_k+s.$$
But $T_t=(T_{\delta_k})^n\circ T_s$.  This implies that if
$\langle x,e_k^*\rangle=0$, then
$$\langle T_t \, x,e_k^*\rangle=0$$
for all $t>0$.  So for any $k \in \Bbb N$, there is a function
$\gamma_k(t)$ such that $$T_t\,e_k=\gamma_k(t)\, e_k,$$
But $T_t$ is an isometry.  We must have $|\gamma_k(t)|=1$ for all
$t>0$.  By the $C_0$ semigroup property of $(T_t)_{t>0}$, it is
easy to see that $\gamma_k(t)=e^{i\omega_k t}$ for some $\omega_k
\in \Bbb R$.  The proof is complete. 
\end{proof}

\end{document}